\documentclass[twocolumn,aps,nofootinbib,secnumarabic,nofootinbib]{revtex4}
\usepackage{amssymb,amsmath2000,amsthm2000,times,mathptmx}
\DeclareMathOperator{\Vol}{Vol\ }
\newcommand{\st}{{\bigm|}}
\newcommand{\dual}{\sp\circ}
\newcommand{\R}{\mathbb{R}}
\newtheorem{theorem}{Theorem}
\newtheorem{corollary}[theorem]{Corollary}

\begin{document}

\title{A low-technology estimate in convex geometry}
\author{Greg Kuperberg}
\thanks{Supported by a Sloan Foundation Graduate Fellowship in Mathematics}
\affiliation{Department of Mathematics, University of Chicago,
    Chicago, IL 60637}
\email[Current email:]{greg@math.ucdavis.edu}

\begin{abstract}
Let $K$ be an $n$-dimensional symmetric convex body with $n \ge 4$ and let
$K\dual$ be its polar body.  We present an elementary proof of the fact that
$$(\Vol K)(\Vol K\dual)\ge \frac{b_n^2}{(\log_2 n)^n},$$
where $b_n$ is the volume of the Euclidean ball of radius 1. The
inequality is asymptotically weaker than the estimate of Bourgain and
Milman, which replaces the $\log_2 n$ by a constant. However, there is
no known elementary proof of the Bourgain-Milman theorem\footnote{The
abstract is adapted from the Math Review by Keith Ball, MR 93h:52010.}.
\end{abstract}

\maketitle

Let $V$ be a finite-dimensional vector space over $\R$ with a
volume element and let $V\sp *$ denote the dual vector space with the
dual volume element.  A \emph{convex body} is a compact convex set with
nonempty interior.  A convex set is \emph{symmetric} if it is invariant
under $x \mapsto -x$.  We define a \emph{ball} to be a symmetric convex
body.  We define $K\dual$, the \emph{dual} of a ball $K \subset V$, by
$$K\dual = \{y \in V\sp * \st y(K) \subset [-1,1]\}.$$
A ball $K$ is an
\emph{ellipsoid} if it is a set of the form $\{x \st \langle x,x
\rangle_K \le 1\}$ for some positive-definite inner product
$\langle\cdot,\cdot\rangle_K$ on $V$.

In this paper we will present a low-technology proof of the following
estimate:

\begin{theorem}
Let $K$ be a symmetric convex body in an $n$-dimensional space $V$
and suppose that there are two ellipsoids $E_1$ and $E_2$
such that $E_1 \subseteq K \subseteq E_2$
and $(\Vol E_2)/(\Vol E_1) = r^n$ with $r \ge 2$. Then
$$\frac{(\Vol K)(\Vol K\dual)}{(\Vol B)(\Vol B\dual)}\ge (2\log_2 r)\sp {-n},$$
where $B$ is an ellipsoid.
\end{theorem}

If $K$ is an arbitrary convex body of dimension $n$, then the largest-volume
ellipsoid $J \subseteq K$, which is called the John ellipsoid, satisfies
$K \subseteq \sqrt{n}K$.  (Proof: If $x \notin \sqrt{n}J$ but $x \in K$, then
$J$ is not the largest ellipsoid in the convex hull of $J\cup \{x,-x\}$.)
It follows that a corollary.

\begin{corollary}  For symmetric convex body $K$ of dimension $n\ge 4$,
$$\frac{(\Vol K)(\Vol K\dual)}{(\Vol B)(\Vol B\dual)}\ge
(\log_2 n)\sp {-n}. $$
\end{corollary}

It is not surprising that this estimate is asymptotically inferior to a
high-technology estimate due to Bourgain and Milman \cite{BM:volume} (see also
Pisier \cite{Pisier:volume}) which says that, for some fixed constant $C$
independent of $n$ and $K$,
$$\frac{(\Vol K)(\Vol K\dual)}{(\Vol B)(\Vol B\dual)}\ge C^{-n}.$$
These estimates can be considered a partial inverse of Santal\'o's inequality
\cite{Santalo:cuerpos}, which states that:
$$(\Vol K)(\Vol K\dual)\le(\Vol B)(\Vol B\dual).$$
There is a nice proof of Santal\'o's inequality due to Saint-Raymond
\cite{Saint-Raymond:volume}.

We begin with some notation which will be used in the proof of the theorem. If
$X$ and $Y$ are two vector spaces, let $P_{X,Y}$ denote the projection from $X
\times Y$ to $Y$ and interpret $X$ and $Y$ as also being the subsets $X \times
\{0\}$ and $\{0\} \times Y$ of $X \times Y$.  If $K$ is a symmetric convex
body, we define the norm $||\cdot||_K$ by setting $||x||_K$ to be the least
positive number $t$ such that $x/t \in K$; in other words, $K$ is the unit ball
of $||\cdot||_K$.  If $A$ and $B$ are two symmetric convex sets in the same
vector space and $p \ge 1$, let
$$A +_p B \stackrel{\mathrm{def}}{=} \{sa+tb \st a \in A, b \in B,
    \mbox{\ and\ } |s|\sp p + |t|\sp p \ \le 1\},$$
and if $A$ and $B$ are convex bodies, let $A \cap_p B$ be the convex body $C$
such that
$$||x||_C\sp p = ||x||_A\sp p + ||x||_B\sp p.$$
(These
definitions are obviously related to the $\ell_p$ norms and have the
usual interpretation when $p = \infty$.) If $A$ is a symmetric convex
set in $X$ and $B$ is a symmetric convex body in $Y$, let $A \times_p B$
denote $A +_p B \subset X \times Y$.  Thus, $+_\infty$, $\cap_\infty$,
and $\times_\infty$ coincide with the usual operations of $+$, $\cap$,
and $\times$ for sets, and $A +_1 B$ is the
convex hull of $A$ and $B$.  Note that the result of any of these
operations is always a symmetric convex body.  Finally, a standard computation
shows that, if $A$ is $n$-dimensional and $B$ is $k$-dimensional,
$$\Vol A \times_p B = \frac{(\Vol A)(\Vol B)}{\binom{(n+k)/p}{n/p}},$$
where a fractional binomial coefficient is interpreted by the factorial
formula, i.e., $x! = \Gamma(x+1)$.

\begin{proof}[Proof of theorem] The result is clearly true if $2 \le r \le
4$, because in this case $r \le 2\log_2 r$, and the volume ratio is at least
$r^{-n}$ because $E\dual_2 \subseteq K\dual$.  Otherwise, let $F$ be the
unique ellipsoid such that if we identify $V$ with $V\sp *$ by the
inner product $\langle\cdot,\cdot\rangle_{F}$, then $E_1 = E_2\dual$.
We will maintain this identification between $V$ and $V\sp *$ for the
rest of the proof, and we can assume to avoid confusion that the volume
elements on $V$ and $V\sp *$ are equal.

Consider the convex body $S(K \times_2 K\dual) \subset V \times V\sp
*$, where $S$ is the linear operator given by $S(x,y) = (x,x+y)$.
Observe that
$$V \cap S(K \times_2 K\dual) = K \cap_2 K\dual$$
and that
$$P_{V,V^*}(S(K \times_2 K\dual)) = K +_2 K\dual.$$
Thus:
\begin{align*} 
(\Vol K)(\Vol K\dual) &= \binom{n}{n/2} (\Vol S(K \times_2 K\dual))\\
    &> \frac{\binom{n}{n/2}}{\binom{2n}{n}}(\Vol K\cap_2 K\dual)
    (\Vol K+_2 K\dual) \\ 
    &> 2^{-n} (\Vol K\cap_2 K\dual)(\Vol (K\cap_2 K\dual)\dual).  
\end{align*}
The first inequality follows from an estimate of Rogers and Shepard:
If $C$ is a symmetric convex body in $X \times Y$, where $X$ and $Y$ are
vector spaces, then $C$ is at least as big as $(C \cap X) \times_1
P_{X,Y}(C)$.  (Proof:  For all $x \in P_{X,Y}(C)$, $C \cap (x+X)$
contains a translate of $(C\cap X)(1-||x||_{P_{X,Y}(C)})$.)

Finally, observe that
$$\frac{\Vol F}{\Vol E_1} =\sqrt{\frac{\Vol E_2}{\Vol E_1}}$$
and that
$$\frac1{\sqrt{2}}F \supseteq K \cap_2 K\dual \supseteq E_1 \cap_2 E_1
    = \frac1{\sqrt{2}}E_1.$$
The first inclusion follows from the observation that
$$||x||_F\sp 2 = \langle x,x\rangle_F < ||x||_K ||x||_{K\dual},$$
which implies that
$$||x||_K\sp 2 + ||x||_{K\dual}\sp 2 \ge 2||x||_F\sp 2.$$
The theorem follows by induction.
\end{proof}

\acknowledgments

I would like to thank the Institut des Hautes \'Etudes Scientifiques
for their hospitality during my stay there.  I would also like to thank
Sean Bates for his encouragement and interest in this work.


\providecommand{\bysame}{\leavevmode\hbox to3em{\hrulefill}\thinspace}

\end{document}